\documentclass[12pt]{amsart}
\usepackage[cp1251]{inputenc}
\usepackage[english,russian]{babel}
\usepackage{amsmath}
\usepackage{amssymb}
\usepackage{amsfonts}
\title{The factorization technique for difference equations}

\author{Alexander V. Vasilyev, Vladimir B. Vasilyev }

\date{}

\begin{document}
\renewcommand{\refname}{References}
\maketitle

{\bf Abstract.} We study  multidimensional difference equations with a continual variable in the Sobolev--Slobodetskii spaces. Using ideas and methods of the theory of boundary value problems for elliptic pseudo differential equations we suggest to consider certain boundary value problems for such difference equations. Special boundary conditions permit to prove unique solvability for these boundary value problems in appropriate Sobolev--Slobodetskii spaces.



\section{Introduction}

We consider a general difference equation of the type
\begin{equation}\label{1}
\sum\limits_{|k|=0}^{\infty}a_k(x)u(x+\alpha_k)=v(x),~~~x\in D,
\end{equation}
where $D\subset\mathbb R^m$ is a canonical domain like $\mathbb R^m, \mathbb R^m_{\pm}=\{x\in\mathbb R^m: x=(x_1,\cdots,x_m), \pm x_m>0\}, C^a_+=\{x\in\mathbb R^m: x_m>a|x'|, x'=(x_1,\cdots,x_{m-1}), a>0\}, k$ is a multi-index, $|k|=k_1+\cdots+k_m, \{\alpha_k\}\subset D$. Equations of a such type have a long history \cite{T,MT,J} and in general there is no an algorithm for solving the equation \eqref{1}. If so then any assertion on a solvability of such equations is very important and required. One can add that equations \eqref{1} appear in very distinct branches of a science like mathematical biology, technical problems, etc. Also such equations have arisen in studies of the second author \cite{V1,V2} related to boundary value problems in a plane corner. One-dimensional case for such equations was considered in \cite{VV4}.

Here we'll start from the equation
\begin{equation}\label{2}
\sum\limits_{|k|=0}^{\infty}a_ku(x+\alpha_k)=v(x),~~~x\in\mathbb R^m_+,
\end{equation}
with constant coefficients because further we'll try to use a local principle \cite{MP} to obtain some results on Fredholm properties of the general equation \eqref{1}. We use methods of the theory of boundary value problems for elliptic pseudo differential equations \cite{E,V}.
For our case of a half-space these methods are based on the theory of one-dimensional singular integral equations and classical Riemann boundary value problem \cite{G,M,GK}.

\section{Spaces, operators and symbols}

\subsection{Spaces}

Let $S(\mathbb R^m)$ be the Schwartz class of infinitely differentiable rapidly decreasing at infinity functions,  $S'(\mathbb R^m)$ be the space of distributions over the space  $S(\mathbb R^m)$. If $u\in S(\mathbb R^m)$ then its Fourier transform is defined by the formula
\[
\tilde u(\xi)=\int\limits_{\mathbb R^m}e^{-ix\cdot\xi}u(x)dx.
\]

{\bf Definition 1.}
{\it A Sobolev--Slobodetskii space $H^s(\mathbb R^m), s\in\mathbb R,$ consists of functions (distributions)
with a finite norm
\[
||u||_s=\left(~\int\limits_{\mathbb R^m}\tilde u(\xi)(1+|\xi|)^{2s}d\xi\right)^{1/2}.
\]
}

Let's note $H^0(\mathbb R^m)=L_2(\mathbb R^m)$.

The space  $S(\mathbb R^m)$ is a dense subspace in the  $H^s(\mathbb R^m)$ \cite{E}. The space $H^s(\mathbb R^m_+)$ consists of functions from the space $H^s(\mathbb R^m)$ which support belongs to $\overline{\mathbb R^m_+}$ with induced norm. Also we need the space $H_0^s(\mathbb R^m_+)$ which consists of distributions from $S'(\mathbb R^m_+)$ admitting a continuation in the whole space  $H^s(\mathbb R^m)$. A norm in the space $H_0^s(\mathbb R^m_+)$ is defined by the formula
\[
||u||^+_s=\inf ||lu||_s,
\]
where {\it infimum} is taken from all continuations $l$.

\subsection{Operators}

Here we consider difference operators with constant coefficients only of the type
\begin{equation}\label{3}
{\mathcal D}: u(x)\longmapsto \sum\limits_{|k|=0}^{\infty}a_ku(x+\alpha_k),
\end{equation}
where $\{a_k\},\{\alpha_k\}$ are given sequences in $\mathbb R^m$, and
\begin{equation}\label{4}
 \sum\limits_{|k|=0}^{\infty}|a_k|<+\infty.
 \end{equation}

{\bf Definition 2.}
{\it An operator $\mathcal D$ of the type \eqref{3} with coefficients $a_k$ satisfying \eqref{4} is called difference operator with constant coefficients.
}

{\bf Lemma 1.}
{\it Every operator ${\mathcal D}: H^s(\mathbb R^m)\rightarrow H^s(\mathbb R^m)$ with constant coefficients is a linear bounded operator $\forall~s\in\mathbb R$.
}

\subsection{Symbols}

{\bf Definition 3.}
{\it The function
\begin{equation}\label{5}
\sigma_{\mathcal D}(\xi)= \sum\limits_{|k|=0}^{\infty}a_ke^{-i\alpha_k\cdot\xi}
\end{equation}
is called a symbol of the operator $\mathcal D$. The symbol $\sigma_{\mathcal D}(\xi)$ is called an {\bf elliptic symbol} if $\sigma_{\mathcal D}(\xi)\neq 0, \forall~\xi\in\mathbb R^m$.
}

Evidently under condition \eqref{4} $\sigma_{\mathcal D}\in L_{\infty}(\mathbb R^m)$, but everywhere below we suppose that $\sigma_{\mathcal D}\in C(\dot{\mathbb R}^m)$ taking into account that $\dot{\mathbb R}^m$ is a compactification of $\mathbb R^m$.

\section{Equations and factorization}

\subsection{Equations}

We are interested in studying solvability of the equation \eqref{2}. It can be written in the operator form
\begin{equation}\label{6}
({\mathcal D}u)(x)=v(x),~~~x\in\mathbb R^m_+,
\end{equation}
assuming that $v$ is a given function in $\mathbb R^m_+$, $v\in H_0^s(\mathbb R^m_+)$, the unknown function $u$ is defined in $\mathbb R^m_+$, $u\in H^s(\mathbb R^m_+)$, and $\{\alpha_k\}\subset\mathbb R^m_+$.

By notation, $u_+(x)=u(x)$, $lv$ is an arbitrary continuation of $v$ on $\mathbb R^m_+$. Then we put
\[
u_-(x)=(lv)(x)-({\mathcal D}u_+)(x),
\]
and see that $u_-(x)=0, \forall~x\in\mathbb R^m_+,$ to explain this notation. Further we rewrite the last equation
\[
({\mathcal D}u_+)(x)+u_-(x)=(lv)(x)
\]
and apply the Fourier transform
\begin{equation}\label{7}
\sigma_{\mathcal D}(\xi)\tilde u_+(\xi)+\tilde u_-(\xi)=\tilde{lv}(\xi).
\end{equation}

To solve the equation \eqref{7} with an elliptic symbol $\sigma_{\mathcal D}(\xi)$ we need to introduce a concept of a factorization. Everywhere below we write $\sigma(\xi)$ instead of $\sigma_{\mathcal D}(\xi)$ for a brevity.

\subsection{Factorization}

Let's denote $\xi=(\xi',\xi_m), \xi'=(\xi_1,\cdots,\xi_{m-1})$.

{\bf Definition 4.}
{\it Let $\sigma(\xi)$ be an elliptic symbol. Factorization of elliptic symbol  $\sigma(\xi)$ is called its representation in the form
\[
\sigma(\xi)= \sigma_+(\xi)\sigma_-(\xi),
\]
where factors  $\sigma_{\pm}(\xi)$ admit an analytical continuation in upper and lower complex planes $\mathbb C_{\pm}$ on the last variable $\xi_m$ for almost all $\xi'\in\mathbb R^{m-1}$ and $\sigma_{\pm}(\xi)\in  L_{\infty}(\mathbb R^m)$.
}

{\bf Definition 5.}
{\it Index of factorization for the elliptic symbol $\sigma(\xi)$ is called an integer
\[
\ae=\frac{1}{2\pi}\int\limits_{-\infty}^{+\infty}d\arg\sigma(\cdot,\xi_m).
\]
}

{\bf Remark 1.}
{\it The index $\ae$ is not really depended on $\xi'$ because it is homotopic invariant.
}

{\bf Remark 2.}
{\it It is principal fact the index of factorization does not correlate with an order of operator. For our case the order of the operator $\mathcal D$ is zero in a sense of G. Eskin's book \cite{E} but the index may be an arbitrary integer. It is essential the index is a {\bf topological barrier} for a solvability.
}

{\bf Proposition 1.}
{\it If $\ae=0$ then for any elliptic symbol $\sigma(\xi)$ a factorization
\[
\sigma(\xi)=\sigma_+(\xi)\sigma_-(\xi)
\]
exists, and it is unique up to a constant.
}

This is classical result, see details in \cite{G,M,GK,E}.

\section{Solvability and boundary value problems}

\subsection{Solvability}

Everywhere below we'll denote $\widetilde H(D)$ the Fourier image of the space $H(D)$.

{\bf Theorem 1.}
{\it If $|s|<1/2,~\ae=0$ then the equation \eqref{6} has a unique solution $u\in H^s(\mathbb R^m_+)$ for arbitrary right hand side $v\in H^s_0(\mathbb R^m_+)$.
}

{\bf Proof} is a very simple. It is based on properties of the Hilbert transform
\[
(H_{\xi'}u)(\xi',\xi_m)=\frac{1}{\pi i}v.p.\int\limits_{-\infty}^{+\infty}\frac{u(\xi',\eta_m)d\eta_m}{\xi_m-\eta_m}
\]
which is a linear bounded operator $H^s(\mathbb R^m)\rightarrow H^s(\mathbb R^m)$ for $|s|<1/2$ \cite{E}. This operator generates two projectors on some spaces consisting of boundary values of analytical functions in $\mathbb C_{\pm}$ on the last variable $\xi_m$ \cite{G,M,GK,E}
\[
\Pi_{\pm}=1/2(I\pm H_{\xi'}),
\]
so that the representation
\[
f=f_++f_-\equiv\Pi_+f+\Pi_-f
\]
is unique for arbitrary $f\in H^s(\mathbb R^m), |s|<1/2$. Further after factorization we write the equality \eqref{7} in the form
\[
\sigma_+(\xi)\tilde u_+(\xi)+\sigma^{-1}_-(\xi)\tilde u_-(\xi)=\sigma^{-1}_-(\xi)\tilde{lv}(\xi),
\]
and else
\[
\sigma_+(\xi)\tilde u_+(\xi)-(\Pi_+(\sigma^{-1}_-\cdot\tilde{lv}))(\xi)=
\]
\[
(\Pi_-(\sigma^{-1}_-\cdot\tilde{lv}))(\xi)-\sigma^{-1}_-(\xi)\tilde u_-(\xi).
\]

So the left hand side belongs to the space $\widetilde H^{s}(\mathbb R^m_+)$ and the left hand side belongs to the space $\widetilde H^{s}(\mathbb R^m_-)$, and these should be zero. Hence
\[
\tilde u_+(\xi)=\sigma^{-1}_+(\xi)(\Pi_+(\sigma^{-1}_-\cdot\tilde{lv}))(\xi).
\]

It completes the proof. $\triangle$

\subsection{General solution}

Let $\ae\in\mathbb Z$. First we introduce a function
\[
\omega(\xi',\xi_m)=\left(\frac{\xi_m-i|\xi'|-i}{\xi_m+i|\xi'|+i}\right)^{\ae},
\]
which belongs to $C(\dot{\mathbb R}^m)$.

Evidently the functions $z\pm i|\xi'|$ for fixed $\xi'\in\mathbb R^{m-1}$ are analytical functions in complex half-planes $\mathbb C_{\pm}$. Moreover
\[
\frac{1}{2\pi}\int\limits_{-\infty}^{+\infty}d\arg\frac{\xi_m-i|\xi'|-i}{\xi_m+i|\xi'|+i}=1.
\]

According to the index property \cite{G,M,GK,E} a function
\[
\omega^{-1}(\xi',\xi_m)\sigma(\xi',\xi_m)
\]
has a vanishing index, and it can be factorized
\[
\omega^{-1}(\xi',\xi_m)\sigma(\xi',\xi_m)=\sigma_+(\xi',\xi_m)\sigma_-(\xi',\xi_m),
 \]
so we have
\[
\sigma(\xi',\xi_m)=\omega(\xi',\xi_m)\sigma_+(\xi',\xi_m)\sigma_-(\xi',\xi_m),
 \]
 where
 \[
 \sigma_{\pm}(\xi',\xi_m)=\exp({\Psi^{\pm}(\xi',\xi_m)}),~~~{\Psi^{\pm}(\xi',\xi_m)}=
\]
 \[
 \frac{1}{2\pi i}\lim_{\tau\to 0+}\int\limits_{-\infty}^{+\infty}\frac{\ln(\omega^{-1}\sigma)(\xi,\eta_m)d\eta_m}{\xi_m\pm i\tau -\eta_m}.
 \]

 Now the equation \eqref{7} we rewrite in the form
\[
(\xi_m+i|\xi'|+i)^{-\ae}\sigma_+(\xi)\tilde u_+(\xi)+(\xi_m-i|\xi'|-i)^{-\ae}\sigma_-^{-1}(\xi)\tilde u_-(\xi)=
\]
\begin{equation}\label{11}
(\xi_m-i|\xi'|-i)^{-\ae}\sigma_-^{-1}(\xi)\tilde{lv}(\xi).
\end{equation}

Let's note the right hand side of the equation \eqref{11} belongs to the space $\widetilde H^{s+\ae}(\mathbb R^m)$.

If $|\ae+s|<1/2$ we go to the sec. 4.1.

\subsubsection{Positive case}

If $ s+\ae>1/2$ we choose a minimal $n\in\mathbb N$ so that $0<s+\ae-n<1/2$.  Further we use a decomposition formula for operators $\Pi_{\pm}$ \cite{E} for $\tilde f\in\widetilde H^{s+\ae}(\mathbb R^m)$
\begin{equation}\label{12}
\Pi_{\pm}\tilde f=\sum\limits_{k=1}^n\frac{\Pi'\Lambda^{k-1}_{\pm}\tilde f}{\Lambda^k_{\pm}}+\frac{1}{\Lambda^n_{\pm}}\Pi_{\pm}\Lambda^n_{\pm}\tilde f,
\end{equation}
where
\[
\Lambda_{\pm}(\xi',\xi_m)=\xi_m\pm|\xi'|\pm i,~~~ (\Pi'\tilde f)(\xi')=\frac{1}{2\pi}\int\limits_{-\infty}^{+\infty}\tilde f(\xi',\xi_m)d\xi_m.
\]

We rewrite the equation \eqref{11}
\[
\sigma_+(\xi)\tilde w_+(\xi)+\sigma_-^{-1}(\xi)\tilde w_-(\xi)=\tilde h(\xi),
\]
where $\tilde w_{\pm}(\xi)=(\xi_m\pm i|\xi'|\pm i)^{-\ae}\tilde u_{\pm}(\xi), \tilde h(\xi)=(\xi_m-i|\xi'|-i)^{-\ae}\sigma_-^{-1}(\xi)\tilde{lv}(\xi).$

Obviously $\tilde w_{\pm}\in\widetilde H^{s+\ae}(\mathbb R^m_{\pm}), \tilde h\in\widetilde H^{s+\ae}(\mathbb R^m)$. We set $s+\ae-n=\alpha, 0<\alpha<1/2$.Since $s+\ae=n+\alpha>\alpha$ then $\tilde h\in\widetilde H^{s+\ae}(\mathbb R^m)\Longrightarrow h\in\widetilde H^{\alpha}(\mathbb R^m)$. According to the theorem 4.1 we have a solution of the last equation $\tilde w_+\in\widetilde H^{\alpha}(\mathbb R^m_+)$ in the form
\[
\tilde w_+(\xi)=\sigma^{-1}_+(\xi)(\Pi_+\tilde h)(\xi).
\]

Thus
\[
\tilde u_+(\xi)=(\xi_m+i|\xi'|+i)^{\ae}\sigma^{-1}_+(\xi)(\Pi_+\tilde h)(\xi),
\]
so that $\tilde u_+\in\widetilde H^{\alpha-\ae}(\mathbb R^m_+)$. Now we apply the formula \eqref{12} to the expression $\Pi_+\tilde h$
and obtain the following representation
\begin{equation}\label{13}
\tilde u_+(\xi)=\sum\limits_{k=1}^n\frac{\tilde c_k(\xi')}{\sigma_+(\xi)\Lambda^{k-\ae}_{+}(\xi',\xi_m)}+\frac{1}{\sigma_+(\xi)\Lambda^{n-\ae}_{+}(\xi',\xi_m)}(\Pi_{\pm}\Lambda^{n}_{+}\tilde h)(\xi',\xi_m),
\end{equation}
where $\tilde c_k=(\Pi'\Lambda^{k-1}_{+})\tilde h$. It is not hard concluding $\tilde c_k\in\widetilde H^{s_k}(\mathbb R^{m-1}), s_k=s+\ae-k+1/2$. So we have the following

{\bf Proposition 2.}
{\it If $s+\ae>1/2$ then for the solution of the equation \eqref{6} the representation \eqref{13} is valid.
}

{\it Note.} {\footnotesize One can prove that the functions $\tilde c_k\in\widetilde H^{s_k}(\mathbb R^{m-1}), s_k=s+\ae-k+1/2,$ are defined uniquely.
}

\subsubsection{Negative case}

If $s+\ae<-1/2$ we choose a polynomial $Q_n(\xi)$ without real zeroes so that $-1/2<s+\ae+n<0$, and use the equality
\[
\sigma_+(\xi)\tilde w_+(\xi)+\sigma_-^{-1}(\xi)\tilde w_-(\xi)=\tilde h(\xi)
\]
from sec. 4.2.1 once again. Since $\tilde h\in\widetilde H^{s+\ae}(\mathbb R^m)$ we represent
\[
\tilde h=Q\Pi_+(Q^{-1}\tilde h)+Q\Pi_-(Q^{-1}\tilde h)
\]
because $Q^{-1}\tilde h\in\widetilde H^{s+\ae+n}(\mathbb R^m)$. Further we work with the equality
\[
\sigma_+(\xi)\tilde w_+(\xi)+\sigma_-^{-1}(\xi)\tilde w_-(\xi)=Q\Pi_+(Q^{-1}\tilde h)+Q\Pi_-(Q^{-1}\tilde h)
\]
or in other words
\[
\sigma_+(\xi)\tilde w_+(\xi)-Q\Pi_+(Q^{-1}\tilde h)=Q\Pi_-(Q^{-1}\tilde h)-\sigma_-^{-1}(\xi)\tilde w_-(\xi)
\]

So the left hand side belongs to the space $\widetilde H^{s+\ae}(\mathbb R^m_+)$ and the left hand side belongs to the space $\widetilde H^{s+\ae}(\mathbb R^m_-)$ so it is distribution supported on $\mathbb R^{m-1}$. Its general form in Fourier images is \cite{E}
\[
\sum\limits_{j=1}^{n}\tilde c_j(\xi')\xi_m^{j-1}.
\]

Thus we have the formula $(\tilde g_+=\Pi_+(Q^{-1}\tilde h)$
\[
(\xi_m+i|\xi'|+i)^{-\ae}\sigma_+(\xi)\tilde u_+(\xi)-Q_n(\xi)g_+(\xi)=\sum\limits_{j=1}^{n}\tilde c_j(\xi')\xi_m^{j-1}.
\]
and a lot of solutions
\[
\tilde u_+(\xi)=(\xi_m+i|\xi'|+i)^{\ae}\sigma^{-1}_+(\xi)Q_n(\xi)g_+(\xi)+
\]
\[
(\xi_m+i|\xi'|+i)^{\ae}\sigma^{-1}_+(\xi)\sum\limits_{j=1}^{n}\tilde c_j(\xi')\xi_m^{j-1}.
\]

It is left to verify that functions $\tilde C_j(\xi)=(\xi_m+i|\xi'|+i)^{\ae}\sigma^{-1}_+(\xi)\tilde c_j(\xi')\xi_m^j$ belong to $\widetilde H^s(\mathbb R^m)$. We have
\[
||C_j||_s^2=\int\limits_{\mathbb R^m}|\tilde c_j(\xi')|^2|\xi_m+i|\xi'|+i|^{2\ae}|\sigma^{-2}_+(\xi)||\xi_m|^{2j}(1+|\xi|)^{2s}d\xi,
\]
and passing to repeated integral we first calculate
\[
\int\limits_{-\infty}^{+\infty}|\xi_m+i|\xi'|+i|^{2\ae}|\xi_m|^{2j}(1+|\xi|)^{2s}d\xi_m,
\]
which exists only if $\ae+j+s<-1/2$. Hence we obtain after integration $\tilde C_j\in H^{\ae+j+s+1/2}(\mathbb R^{m-1})$.

Thus we have proved the following

{\bf Theorem 2.}
{\it If $s+\ae<-1/2$, then the equation \eqref{6} has many solutions in the space $H^s(\mathbb R^m_+)$, and formula for a general solution in Fourier image
\[
  \tilde u_+(\xi',\xi_m)=(\xi_m+i|\xi'|+i)^{\ae}\sigma^{-1}_+(\xi)Q_n(\xi)\tilde g_+(\xi',\xi_m) +
\]
\[
(\xi_m+i|\xi'|+i)^{\ae} \sigma^{-1}_+(\xi',\xi_m) \sum\limits_{k=0}^{\ae-1}c_k(\xi')\xi_m^k
\]
holds, where  $c_k\in H^{s_k}(\mathbb R^{m-1}), s_k=-\ae+k+1/2, k=0,\cdots,\ae-1,$ are arbitrary functions.
}

{\bf Corollary.}
{\it If  under assumptions of the theorem 4.3 $v\equiv 0$ then a general solution of the equation
\begin{equation}\label{14}
({\mathcal D}u)(x)=0,~~~~x\in\mathbb R^m_+
 \end{equation}
 has the form
 \begin{equation}\label{15}
 \tilde u_+(\xi)= (\xi_m+i|\xi'|+i)^{\ae} \sigma^{-1}_+(\xi',\xi_m) \sum\limits_{k=1}^{n}\tilde c_k(\xi')\xi_m^{k-1}.
 \end{equation}
}

\subsection{Boundary conditions}

For a brevity we consider a homogeneous equation using the corollary 4.4. We need some additional conditions to uniquely determine the functions $\tilde c_k, k=1,\cdots,n$. It is interesting fact that we can't use the same conditions for positive and negative $\ae$. Moreover the boundary operators in a certain sense are determined by the formula for a general solution. We consider below very simple boundary operators. Usually such operators are traces of some pseudo differential operators on the hyper-plane $x_m=0$. But it is possible not for all cases.
\subsubsection{Positive case}

Let us assume we know  the values of $\tilde u_+$ in $n$ distinct hyper-planes from $\mathbb R^m$ of type $\xi_m=p_j$. We denote
$\tilde u_+(\xi',p_j)\equiv\tilde r_j(\xi')$ and obtain from the formula \eqref{15} the following system of linear algebraic equations
\[
\sum\limits_{k=1}^{n}\tilde c_k(\xi')p_j^{k-1}=\tilde r_j(\xi') (p_j+i|\xi'|+i)^{\ae} \sigma^{-1}_+(\xi',p_j),~~~~~~j=1,\cdots,n.
\]

Obviously the system is uniquely solvable because its matrix has the Vandermonde determinant. To formulate a corresponding boundary value problem we need some preliminaries.

We take the following boundary conditions
\begin{equation}\label{16}
\int\limits_{-\infty}^{+\infty}u_+(x',x_m)e^{-ip_jx_m}dx_m=r_j(x'),~~~j=1,\cdots,n.
\end{equation}

It will mean $\tilde u_+(\xi',p_j)=\tilde r_j(\xi')$. If $u_+\in H^s(\mathbb R^m_+)$ then $r_j\in H^{s-1/2}(\mathbb R^m_+)$ \cite{E}. So we have the following

{\bf Theorem 3.}
{\it Let $r_j\in H^{s-1/2}(\mathbb R^{m-1}), j=1,\cdots,n$. Then the boundary value problem \eqref{14},\eqref{16} has a unique solution in the space $H^s(\mathbb R^m_+)$.
}

{\it Note.} {\footnotesize One can consider a linear combination of the conditions \eqref{16} and require non-vanishing the associated determinant.
}

\subsubsection{Negative case}

This case admits integration for the right hand side of the formula \eqref{15} thus we take boundary conditions in the standard form
\begin{equation}\label{17}
(A_ju_+)(x)_{|_{x_m=0}}=r_j(x'),~~~j=1,\cdots,n,
\end{equation}
where $A_j$ are pseudo differential operators with symbols $A_j(\xi',\xi_m)$ satisfying the condition
\[
|A_j(\xi',\xi_m)\sim (1+|\xi'|+|\xi_m|)^{\gamma_j}.
\]

Let's denote
\[
a_{jk}(\xi')=\int\limits_{-\infty}^{+\infty}A_j(\xi',\xi_m)(\xi_m+i|\xi'|+i)^{\ae} \sigma^{-1}_+(\xi',\xi_m) \xi_m^{k-1}d\xi_m.
\]

{\bf Theorem 4.}
{\it Let $\gamma_j+\ae+k<-1, r_j\in H^{s_j}(\mathbb R^{m-1}), s_j=s-\gamma_j-1/2,~\forall j, k=1,\cdots,n,$ and the
\[
\inf_{\xi'\in\mathbb R^{m-1}}|\det(a_{jk}(\xi'))_{j,k=1}^n|>0.
\]
Then the boundary value problem \eqref{14},\eqref{17} has a unique solution in the space $H^s(\mathbb R^m_+)$.
}

\section*{Conclusion}

There are a lot of possibilities to state distinct problems for the equation \eqref{6} adding some additional conditions. Also it seems to be interesting to transfer this approach and results to a discrete case, i.e. for spaces of a discrete variable. This will be discussed elsewhere.


\subsection*{Acknowledgment}
This work was partially supported by Russian Foundation for Basic Research and government of Lipetsk region of Russia, project no. 14-41-03595-a.

\newpage

Alexander V. Vasilyev\\
Studencheskaya 14/1\\
National Research Belgorod State University\\
Belgorod 307008\\
Russia\\
e-mail: alexvassel@gmail.com\\

Vladimir B. Vasilyev\\
Moskovskaya 30\\
Lipetsk State Technical University\\
Lipetsk 398600\\
Russia\\
e-mail: vbv57@inbox.ru

\end{document}